
\input amstex
\documentstyle {amsppt}
\magnification=\magstep1
\define\id{\operatorname{id}}
\define\tr{\operatorname{tr}}\topmatter
\define\ad{\operatorname{ad}}
\title
On von Neumann algebras\\
 which are complemented subspaces of $B(H)$
\endtitle
\author
\\
{\rm{by}}\\
\\
Erik Christensen and Allan M\. Sinclair
\endauthor
\endtopmatter

\document
\subheading{1. Introduction}

A von Neumann algebra $\Cal M$ on a separable Hilbert space $H$ is said to
be approximately finite dimensional if it is the weak closure of an
increasing sequence of finite dimensional $C^\ast$-algebras. These
algebras have played an important role in the study of general von
Neumann algebras, since they have a rich structure and contain a
subgroup of the unitary group which is amenable and still generate
the entire von Neumann algebra. It is fairly easy to prove that an
approximately finite dimensional von Neumann algebra $\Cal M$ is
injective in the sense that there exists a completely positive
projection $\pi$ of $B(H)$ onto $\Cal M$. Such a projection has norm one
and is an $\Cal M$-module mapping. One of Alain Connes major achievements
consists in proving the opposite. In [3] he proved that an injective
factor on a separable Hilbert space is approximately finite
dimensional. It is a natural question then to ask whether the
projection onto $\Cal M$ has to have norm one and be completely positive
in order to prove approximately finite dimensionality.  Bunce and
Paschke [1] and G\. Pisier [11] have some partial results which
indicate that the question may have a positive answer. The results
here are based on [1] and the dependence will be described later.
Pisiers work [11] shows that if a von Neumann
algebra $\Cal M$ is complemented as subspace of $B(H)$ and as a von
Neumann algebra  is isomorphic to the von Neumann algebra tensor
product $\Cal M\overline\otimes \Cal M$ then $\Cal M$
is injective, i.e\. there
exists a norm one projection onto $\Cal M$.

After finishing the first draft of this paper we obtained a
preliminary version of [12] by G\. Pisier where he among other things
proves that if there exists a completely bounded projection from
$B(H)$ onto a von Neumann algebra $\Cal M$ on $H$ then $\Cal M$
is injective.
Since we have obtained the same result but with different methods and
since we think our methods may have some interest in their own, we
have decided to complete the manuscript for publication.
In both papers the completely bounded projection case is the major
one. From here both G\. Pisier and we can obtain various results on
the general case where we just have a bounded projection onto $\Cal M$.

The structure of the paper is as follows: In Lemma 2.1 we present our
key result which we hope may  be usefull in the future studies
of completely bounded operators.  The lemma is based on techniques
developped in [2, 8, 13] and the proof is essentially a fixed point
argument based on a trick due to Sakai [13]. The result is used to
prove that if there exists a completely bounded projection of $B(H)$
onto $\Cal M$ then there exists an $\Cal M$-module projection of $B(H)$ onto
$\Cal M$. We are now exactly in the situation described in [1], and
following Bunce and Paschke we conclude that $\Cal M$ is injective.

The case of a norm continuous projection can be dealt with in a
number of cases, once we have the completely bounded projection case.
If $\Cal M$ is a properly infinite von Neumann algebra, or $\Cal M$
is a finite
von Neumann algebra of type II$_1$,
 which can be decomposed into a tensor product
$\Cal S\overline\otimes\Cal T$ of two von Neumann algebras of type
II$_1$, it is possible to use the results obtained so far to prove
that in these cases, $\Cal M$ is injective if $\Cal M$ is complemented
in $B(H)$. These investigations have two other offsprings. The first,
Proposition 2.3 says that if $x,y$ are bounded operators on a Hilbert
space $H$, $\Cal M$ is a von Neumann algebra on $H$, $x\Cal M
y\subseteq\Cal M$ and $x$ is positive and invertible, then there
exists operators $n_1,n_2$ in $\Cal M$ such that for each $m$ in
$\Cal M$ $xmy=n_1mn_2$. We hope this result may be the first in a
series of result which will help our understanding of the structure
of normal completely bounded operators on a von Neumann algebra. The
other off spring we have got, is that we have now understood Connes
paper [4] and found out, that for a von Neumann algebra $\Cal M$ in
standard position on a Hilbert space $H$ we can equip the space $Y$
consisting of all bounded $\Cal M$-module mappings on $B(H)$, which
vanish on $\Cal M$ with an $\Cal M$-module structure, such that $\Cal
M$ is injective if and only if the derivation $\delta:\Cal M\to Y$
given by $\delta(m)(x)=Jm^\ast Jx-xJm^\ast J$ is inner, (here $J$ is
an antilinear isometric involution on $H$ such that $J\Cal M
J=\Cal M'$, the commutant of $\Cal M$).

\subheading{2. Averaging Completely Bounded Mappings}

The present section contains a fairly general result which we think
may have some interest of its own besides the application made in
Section 3. The result is at least in principle related to Dixmier's
Approximation Theorem [5, III \S5] or [9, 8.3]. The techniques used
are similar to the ones used in [2], which in turn are based on ideas
from [8], [13] by Kadison and Sakai.

\proclaim{Lemma 2.1} Let $\Cal M$ be a von Neumann algebra on a
Hilbert space $H$, and $\goth A$ a
$C^\ast$-algebra on $H$ containing $\Cal M$,
$\pi$ a $\ast$-representation of  $\goth A$ on a Hilbert space $K$
and $S,T$ bounded operators  $S\in B(K,H)$, $T\in B(H,K)$ such that
for elements $a$ in $\goth A$
$$
 S\pi(a)T\in\Cal M\;.
$$
Let $\Cal P$ denote the set of all families $(m_i)_{i\in J}$ of
elements from $\Cal M$ such that $\sum m_i^\ast m_i=I$.
Then for each family $(m_i)$ from $\Cal P$ the sum
$\sum\pi(m_i^\ast) T\,m_i$ converges ultrastrong $-\ast$ in $B(H,K)$.

There exists an operator $\hat T$ in the ultraweakly
closed convex hull of the set $\{\sum\pi(m_i)^\ast
T\,m_i\mid(m_i)\in\Cal P\}$ such that for elements $x$ in $\goth A$
and $m$ in $\Cal M$
 $$S\pi(x)\hat T\in\Cal M\quad\text{ and }\quad S\pi(xm)\hat
T=S\pi(x)\hat T m\;.
 $$
\endproclaim

\demo{Proof} Let $(m_i)_{i\in J}$ be an element of $\Cal P$, then we
can construct a \lq\lq column" operator $v$, which is an isometry of
$H$ into $\bigoplus_J H$ by $v\xi=(m_i\xi)$. Since a
$\ast$-representation is completely bounded the operator
$w=(\pi\otimes id)(v)$ is a contraction operator in $B(K,
\bigoplus_J K)$ given by $w\eta=(\pi(m_i)\eta)$. Let $\tilde T$
denote the operator on $B(\bigoplus_J H,\bigoplus_J K)$
given by $\tilde T(\rho_i)=(T\rho_i)$, then
$\|\tilde T\|=\|T\|$ and it is elementary to check that the sum
$\sum_J\pi(m_i)^\ast T\,m_i$ converges ultrastrong $-\ast$ to
$w^\ast\tilde T v$ in $B(H,K)$. In fact we find that we may define an
operator $\Psi_{(m_i)}:B(H,K)\to B(H,K)$ by
$\Psi_{(m_i)}(x)=w^\ast\tilde x v$ for each element $(m_i)$ in $\Cal
P$. Clearly $\|\Psi_{(m_i)}\|\leq1$ and each $\Psi_{(m_i)}$ is
ultraweak to ultraweak continuous.

Let now $(m_i)_J$ and $(n_l)_L$ be families in $\Cal P$ and let $C$
denote the operator $\sum_L\pi(n_l)^\ast T\,n_l$, then since the sums
involved below are ultraweakly convergent we get $\sum_J\sum_L
m_i^\ast n_l^\ast n_lm_i=I$ and
 $$\Psi_{(m_i)}(C)=\sum_J\pi(m_i)^\ast(\sum_L\pi(n_l)^\ast
Tn_l)m_i=\sum_{J\times L}\pi(n_lm_i)^\ast T\,n_lm_i\;.
 $$

If we now let $\Cal C$ denote the ultraweak closure of the set of
operators $\{\,\sum_J\pi(m_i)^\ast T\,m_i\mid(m_i)\in\Cal P\,\}$,
then the previous arguments show that $\Cal C$ is invariant under all
the operators $\Psi_{(m_i)}$. It is straight forward to check that
$\Cal C$ is convex and by construction for each element $C$ in $\Cal
C$, $\|C\|\leq\| T\|$ so $\Cal C$ is also ultraweakly compact.

By the invariance of $\Cal C$ under $\Psi_{(m_i)}$, $(m_i)\in\Cal P$,
the ultraweak continuity of $\Psi(m_i)$ and the ultraweak compactness
we can conclude that there exists an ultraweakly closed, convex
subset $\Cal C_0$ of $\Cal C$ which is minimal with respect to the
property that $\Psi_{(m_i)}(\Cal C_0)\subseteq\Cal C_0$ for each
$(m_i)$ in $\Cal P$. Let $\hat T\in\Cal C_0$ and let $e$ be a
projection in $\Cal M$ then if $c(e)$ denotes the central support
projection for $e$ in $\Cal M$ we will prove that
 $$\|\hat T e\|=\|\pi(e)\hat T e\|=\|\hat T c(e)\|\;.\tag1
 $$
Since all the operators $\Psi_{(m_i)}$ on $B(H,K)$ are contractions
and $\Cal C_0$ is minimal invariant we must have for each central
projection $g;\|\hat Tg\|=\|\Psi_{(m_i)}(\hat T)g\|$. Choose a family
of partial isometries $(v_j)$ from $\Cal M$ such that $v_jv_j^\ast\leq e$
and $\sum v_j^\ast v_j=c(e)$ then the family
$((I-c(e)),(v_j))$ is in $\Cal P$ so
 $$\align
\|\hat Tc(e)\|&=\|[\pi(I-c(e))\hat T(I-c(e))+\sum\pi(v_j)^\ast\hat T
v_j]c(e)\|\\
&=\max_j\|\pi(v_j)^\ast\hat T v_j\|\\
&\leq\|\pi(e)\hat T e\|\\
&\leq\|\hat T e\|\\
&\leq\|\hat Tc(e)\|\;.\endalign
 $$
 So (1) follows.

Suppose $e$ is a projection in $\Cal M$ and $x$ an operator in $\goth
A$ such that $a=S\pi(x)\pi(I-e)\hat Te\neq0$ then $a$ is in $\Cal M$
so $a^\ast a\in\Cal M$ and we can define a non zero projection $f$ in
$\Cal M$ as the spectral projection for $a^\ast a$ corresponding to
the closed interval $[\frac12\|a\|^2,\|a\|^2]$. Let $\varepsilon$
denote the positive number
 $$\varepsilon=\|a\|^2(4\|S\pi(x)\|^2\|\hat T f\|^2)^{-1}
 $$
 and choose a unit vector $\xi$ in $fH$ such that
 $$\|\pi(f)\hat T f\xi\|^2\geq(1-\varepsilon)\|\hat T f\|^2\;,
 $$
 then:
 $$\align
(1-\varepsilon)\|\hat T f\|^2&\leq\|\pi(f)\hat T f\xi\|^2\\
&\leq\|\pi(e)\hat T f\xi\|^2\\
&=\|\hat T f\xi\|^2-\|(I-\pi(e))\hat T f\xi\|^2\\
&\leq\|\hat T f\|^2-\|(I-\pi(e))\hat T f\xi\|^2\\
&\leq\|\hat T f\|^2-\|S\pi(x)\|^{-2}\|S\pi(x)\pi(I-e)\hat T
ef\xi\|^2\\
&=\|\hat T f\|^2-\|S\pi(x)\|^{-2}\|af\xi\|^2\\
&\leq\|\hat T f\|^2-\|S\pi(x)\|^{-2}\|a\|^2/2\\
&=\|\hat Tf\|^2-2\varepsilon\|\hat T f\|^2\\
&=(1-2\varepsilon)\|\hat T f\|^2\;.\endalign
 $$

We have now obtained a contradiction to the assumption $a\neq0$, so
for each projection $e$ in $\Cal M$ and each $x$ in $\goth A$ we get:
 $$S\pi(x)\pi(I-e)\hat Te=0\;.\tag2
 $$

Let $Q$ denote the orthogonal projection from $K$ onto the closed
linear span of the set $\{\pi(\goth A)S^\ast H\}$. Then $Q$ is in the
commutant $\pi(\goth A)'$ and we get
 $$\forall\,x\in\goth A:\quad S\pi(x)Q\hat T=S\pi(x)\hat T\;.\tag3
 $$

 $$\forall\text{ projections }e\text{ in }\Cal M:\pi(I-e)Q\hat T
e=Q\pi(I-e)\hat T e=0\;.\tag4
 $$
 $$Q\hat T e=\pi(e)Q\hat T e=\pi(e)Q\hat T\;.\tag5
 $$

Since the linear span of projections is norm dense in the von Neumann
algebra $\Cal M$ we get
 $$\forall\, m\in\Cal M:Q\hat T m=\pi(m)Q\hat T\;.\tag6
 $$

A combination of (3) and (6) show that
 $$\forall\, x\in\goth A\;\forall\, m\in\Cal M:\quad S\pi(xm)\hat
T=S\pi(x)\hat Tm\;,\tag7
 $$
 so the lemma follows.$\qquad\square$
\enddemo

\proclaim{Proposition 2.2} Let $\Cal M\subseteq\goth A\subseteq B(H)$
be a pair consisting of a von Neumann algebra $\Cal M$ contained in a
$C^\ast$-algebra $\goth A$ of bounded operators on a Hilbert space
$H$, and let $\Psi$ be a completely bounded operator of $\goth A$
into $\Cal M$. Let $\Cal P$ denote the set of all families
$(m_i)_{i\in J}$ of elements from $\Cal M$ such that $\sum m_i^\ast
m_0=I$, then for each $p=(m_i)_{i\in J}$ there exist completely
bounded operators ${}_p\Psi$ and $\Psi_p$ from $\goth A$ into
$\Cal M$ such that
 $$\multline
\|_p\Psi\|_{cb}\leq\|\Psi\|_{cb}\;,\;\|\psi_p\|_{cb}\leq\|\Psi\|_{cb}\\
{}_p\Psi(a)=\sum_i
m_i^\ast\Psi(m_ia)\;,\;\Psi_p(a)=\sum_i\Psi(am_i^\ast)m_i\;,\endmultline
 $$
 and the sums are ultrastrong $-\ast$ convergent.

Let $L_\Psi$ (resp\. $R_\Psi$) denote the closures in the topology of
pointwise ultraweak convergence of the sets
 $$\{{}_p\Psi\mid p\in\Cal P\,\}\,,\;(\text{resp. }\{\,\Psi_p\mid
p\in\Cal P\,\})\;.
 $$
Then $L_\Psi$ (resp\. $R_\Psi$) contains a completely bounded left
$\Cal M$-module (resp\. right $\Cal M$-module) mapping of $\goth A$
into $\Cal M$.
\endproclaim

\demo{Proof} The proposition follows directly from Lemma 2.1 and the
characterization of completely bounded mappings [10, p\. 105],
[16].$\qquad\square$
\enddemo

The problem of finding criteria on bounded operators $x,y$ on a Hilbert
space $H$ and a von Neumann algebra $\Cal M$ on $H$ such that $x\Cal
M y\subseteq\Cal M$ is only understood to some extent when
$x=y^\ast$.
If $y^\ast\Cal M y\subseteq\Cal M$ and $y$ is invertible, we find
that $y^\ast y\in\Cal M$, so for the polar decomposition $y=vh$ we
get $h$ is a positive invertible in $\Cal M$ and $v$ is a unitary
such that $v^\ast\Cal M v\subseteq\Cal M$. Hence the mapping $m\to
y^\ast my$ has a kind of polar decomposition
$m\overset\Psi\to\rightarrow v^\ast mv\overset\psi\to\rightarrow
hv^\ast mvh$, where $\Psi$ is an \lq\lq inner" completely positive
mapping on $\Cal M$. The Lemma 2.1 above can be used to extend our
information concerning completely bounded operators on a von Neumann
algebra a little bit.

\proclaim{Proposition 2.3} Let $\Cal M$ be a von Neumann algebra on a
Hilbert space  $H$ and let $x,y$ be bounded operators on $H$ such
that $x\Cal M y\subseteq\Cal M$. If $x$ is positive and invertible
then there exist operators $n_1,n_2$ in $\Cal M$ such that for $m$ in
$\Cal M$ $xmy=n_1 mn_2$.
\endproclaim

\demo{Proof} By Lemma 2.1  there exists a positive operator $h$ in
the ultraweakly closed (convex) hull of the set $\{\sum m_i^\ast
xm_i|(m_i)\subseteq\Cal M, \sum m_i^\ast m_i=I\}$ such that for
each $m,n$ in $\Cal M$ $hmy\in\Cal M$ and $hmny=mhny$. Since
$x\geq\delta I$ for some $\delta>0$ we get $\sum m_i^\ast
xm_i\geq\delta I$ so $h\geq\delta I$.

Let $q$ denote the projection in the commutant $\Cal M'$ of $\Cal M$
which is given as the orthogonal projection of $H$ onto the closed
linear span of $\{\Cal M yH\}$, then the relation $mhny=hmny$ implies
that $hq$ is in the commutant $\Cal M'$. Since $h$ is invertible we
get
 $$\align
\forall\,m\in\Cal M: xmy&=xh^{-1}hmqy\\
&=xh^{-1}hqmy\\
&=xh^{-1}mhqy\\
&=xh^{-1}mhy\;.\endalign
 $$
 Now $hy$ is in $\Cal M$ so we define $n_2=hy$.

For any $m$ in $\Cal M$ we have $xh^{-1}mn_2$ is in $\Cal M$ so for
the range projection $R[mn_2]$ of $mn_2$ it follows that
$xh^{-1}R[mn_2]$ is in $\Cal M$. Since the rangeprojection of the sum
of a finite set of projections equals the sup over the same set of
projections we get $xh^{-1}(\underset{i=1}\to{\overset k\to
\vee}R[m_in_2])$ is in $\Cal M$ and consequently for the central support
projection $c$ in $\Cal M$ for $n_2$ we get $xh^{-1}c$
is in $\Cal M$. Let $n_1=xh^{-1}c$ then;
 $$xmy=xh^{-1}mn_2=xh^{-1}mcn_2=xh^{-1}cmn_2=n_1mn_2\;,
 $$
 and the proposition follows.$\qquad\square$
\enddemo

Our next result is a direct application of Lemma 2.1 which will be
extremely usefull in the next section.

\proclaim{Lemma 2.4} Let $\Cal M\subseteq B(H)$ be a von Neumann
algebra and suppose $\rho$ is a completely bounded projection of
$B(H)$ onto  $\Cal M$, then there exists a completely bounded $\Cal
M$-module projection $\gamma$ of $B(H)$ onto $\Cal M$.
\endproclaim

\demo{Proof} Let $\rho(\cdot)=S\pi(\cdot)T$ where $\pi$ is a
$\ast$-represen\-ta\-ti\-on of $\Cal M$ on a Hil\-bert spa\-ce $K$ and
$S^\ast,T\in B(H,K)$. By Lemma 2.1 there exists a $\hat T$ in the
ultraweakly closed convex hull in $B(H,K)$ of the set
$\{\sum_J\pi(m_i)^\ast Tm_i\mid(m_i)_J\in\Cal P\}$ such that
$S\pi(xm)\hat T=S\pi(x)\hat T m$.

Since $\rho$ is a projection we get for each operator
$T_0=\sum_J\pi(m_i)^\ast Tm_i$ from the set and each $m_0$ from $\Cal
M$
 $$\align
S\pi(m_0)T_0&=S(\,\sum_J\pi(m_0m_i^\ast)
Tm_i)\\
&=\sum_J\rho(m_0m_i^\ast)m_i=\sum_Jm_0m_i^\ast m_i=m_0.\endalign
$$

Hence $S\pi(m_0)\hat T=m_0$ and $S\pi(xm_0)\hat T=S\pi(x)\hat Tm_0$.

If we then repeat the procedure with respect to $S$ instead of $T$ we
find an $\hat S$ such that
 $$\multline
\forall\,m\in\Cal M\forall\, x\in B(H)\;.\;\hat S\pi(m)\hat
T=m\,,\;\hat S\pi(mx)\hat T=m\hat S\pi(x)\hat T\;,\\
\hat S\pi(xm)\hat T=\hat S\pi(x)\hat T m\,,\;\hat S\pi(x)\hat
T\in\Cal M\endmultline
 $$

If we then define $\gamma:B(H)\to\Cal M$ by $\gamma(x)=\hat
S\pi(x)\hat T$, we have found a completely bounded $\Cal M$-bimodule
projection of $B(H)$ onto $\Cal M$.$\qquad\square$
\enddemo

 \subheading{3. Characterization of some injective von Neumann
algebras as complemented subspaces of $B(H)$}

Tomiyama studied  von Neumann algebras with \lq\lq extension
property", meaning that we consider a von Neumann algebra on a
Hilbert space $H$ such that there is a projection of norm one from
$B(H)$ onto $\Cal M$, [15]. Among other things he proved that such a
projection is necessarily completely positive. Working in the
category of von Neumann algebras and completely positive mappings the
von Neumann algebras with extension property are the injective
objects. Thanks to Connes fundamental paper [3] it turned out that an
injective factor $\Cal M$ on a separable Hilbert spaces $H$ is the
ultraweak closure of an increasing sequence of finite-dimensional
$C^\ast$-algebras, or otherwise formulated an injective factor on a
separable Hilbert space is approximately finite dimensional. Such
algebras are especially nice because the unitary group $\Cal U$
contains a subgroup $\Cal U_0$ which generates the von Neumann
algebra and is the union of an increasing sequence of compact groups
(i.e\. the unitary groups in the sequence of finite-dimensional
subalgebras).
This subgroup $\Cal U_0$ is clearly amenable and it turns out that for
a general von Neumann algebra $\Cal M$ on a Hilbert space, $\Cal M$
is injective if and only if the unitary group contains an amenable
subgroup which generates $\Cal M$, [6].

Our main result which we are going to present here shows that a von
Neumann algebra $\Cal M$ is injective if there exists a completely
bounded projection onto $\Cal M$ or if there exists a bounded $\Cal
M$-module projection onto $\Cal M$. The argument which yields
injectivity under these conditions comes from [1].

\proclaim{Theorem 3.1} Let $\Cal M$ be a von Neumann algebra
on a Hilbert space $H$. The
following conditions are equivalent

\roster
 \item"i)"  $\Cal M$ is injective
\smallskip
\item"ii)"  There exists a completely bounded projection of $B(H)$
onto $\Cal M$.
\smallskip
 \item"iii)"  There exists a bounded $\Cal M$-module projection of
$B(H)$ onto $\Cal M$.
 \endroster
\endproclaim

\demo{Proof} The implication i) $\Rightarrow$ ii) is already
mentioned since a projection of norm one is completely positive and
hence completely bounded. The result ii) $\Rightarrow$ iii) is proved
in Lemma 2.4. The step iii) $\Rightarrow$ i) is the main result of
the paper [1] by Bunce and Paschke. In the special case where $\Cal
M$ is a finite factor the argument comes directly from Connes paper
[3] in a very transparent way, so we will reproduce it here.
 Suppose
 $\Cal M$ is a
factor of type II$_1$ on a separable Hilbert space $H$ with a cyclic
and separating trace vector $\xi$ and let $\pi$ be a bounded $\Cal
M$-module projection of $B(H)$ onto $\Cal M$. The functional $\psi$
on $B(H)$ given by
 $$\psi(x)=(\pi(x)\xi,\xi)
 $$
 satisfies the property
 $$\align
\forall\,x\in B(H)\;,\;\forall\, m\in\Cal
M:\psi(xm)&=(\pi(xm)\xi,\xi)\\
&=\tr(\pi(x)m)\\
&=\tr(m\pi(x))\\
&=\psi(mx)\;.\endalign
 $$
 Replacing $\psi$ by $\frac12(\psi+\psi^\ast)$ does not affect the
identity above. Since for each unitary $u$ in $\psi$ we have
$\psi(uxu^\ast)=\psi(x)$, the uniqueness of the decomposition of
$\psi$ as a difference of positive functionals $\psi=\psi^+-\psi^-$
such that $\|\psi\|=\|\psi^+\|+\|\psi^-\|$ shows that
$\psi^+(uxu^\ast)=\psi^+(x)$. Since $\psi^+(I)\geq\psi(I)=1$ we
get that the functional $\varphi$ given by
$\varphi=\|\psi^+\|^{-1}\psi^+$ is a hypertrace meaning that
$\varphi\mid\Cal M=\tr$ and for $x$ in $B(H)$ and $m$ in $\Cal M$
$\varphi(xm)=\varphi(mx)$.  By [3, Theorem 5.1,  $7\Rightarrow1$.] We
find that $\Cal M$ is the hyperfinite II$_1$ factor.

If $\Cal M$ is finite but not a factor one can still use the argument
above but one has to refer to Haagerup [7, Lemma 2.2] where he proves
that $\Cal M$ is injective if there exists a family of hypertraces
which separates the central projections in $\Cal M$. The general case
in [1] is done by referring to the decomposition of a type III
algebra as a crossed product of a semifinite algebra and a one
parameter group of automorphisms.$\qquad\square$
\enddemo

\proclaim{Corollary 3.2} If $\Cal M$ is properly infinite or $\Cal M$
is of type II$_1$ and $\Cal M$ is isomorphic to a tensorproduct $\Cal
S\overline\otimes\Cal F$  where $\Cal S$ and $\Cal T$ both are of
type II$_1$, then $\Cal M$ is injective if and only if there exists a
bounded projection of $B(H)$ onto $\Cal M$.
\endproclaim

\demo{Proof} Suppose first $\Cal M$ is properly infinite, let
$\pi:B(H)\to\Cal M$ be a bounded projection onto $\Cal M$ and let
$\Cal S\subseteq\Cal M$ be a von Neumann subalgebra of $\Cal M$
isomorphic to $B(l^2(\Bbb N))$. Finally let $\Cal N$ denote the
relative commutant of $\Cal S$ in $\Cal M$, $\Cal N=\Cal S'\cap\Cal
M$. Now $\Cal M$ is isomorphic to the von Neumann algebra tensor
product $\Cal N\overline\otimes\Cal S$, and we will just use this
description along with the assumptions that $\Cal N$ and $\Cal S$ are
commuting subalgebras of $\Cal M$, generating $\Cal M$. Inside $\Cal
S$ we consider the unitary subgroup $\Cal U_0$ consisting of infinite
unitary matrices which are finite rank perturbations of the identity
I. The group $\Cal U_0$ is clearly the union of an increasing
sequence of compact groups and hence amenable so if for $x$ in $B(H)$
we look at the ultra weakly closed convex hull of the set
 $$\{\, u\pi(u^\ast x)\mid u\in\Cal U_0\,\}\;,
 $$
 we find that we can define an operator $\rho:B(H)\to\Cal M$ such
that $\rho\mid\Cal M=\id\mid \Cal M$
and $\forall\,u\in\Cal U_0$  $\rho(ux)=u\rho(x)$. If we repeat the
arguments, but now to the right we get a bounded projection $\gamma$
of $B(H)$ onto $\Cal M$ such that
 $$\forall\,u\in\Cal U_0\forall\,x\in
B(H):\rho(ux)=u\rho(x)\;,\;\rho(xu)=\rho(x)u\;.
 $$
 Let $(e_{ij})_{i,j\in\Bbb N}$ denote the cannonical matrix
units for $\Cal S$, then the relations above show that
 $$\align
\forall\,x\in B(H):\rho(e_{11}xe_{11})&=e_{11}\rho(x)e_{11}\\
\rho(e_{ii}xe_{jj})&=e_{i1}\rho(e_{1i}xe_{j1})e_{1j}\;.\endalign
 $$
 This means that we can define a projection $\rho_{11}$ of
$B(e_{11}H)$
onto $e_{11}\Cal M\mid e_{11}H$ and that $\rho$ is nothing but
$\rho_{11}\otimes\id$.  In particular $\rho_{11}$ and hence $\rho$ is
completely bounded so we have obtained a completely bounded
projection of $B(H)$ onto $\Cal M$ and $\Cal M$ is then injective by
the theorem.

If $\Cal M$ is finite of type II$_1$ and $\Cal M$ is isomorphic to
the von Neumann algebra tensor product of 2 type II$_1$ von Neumann
algebra $\Cal S$ and $\Cal T$, then both $\Cal S$ and $\Cal T$
contain an injective subfactor $\Cal R_{\Cal S}$ and $\Cal R_{\Cal
T}$ of type II$_1$. We will show that the algebras $\Cal
S\overline\otimes\Cal R_{\Cal T}$ and $\Cal R_{\Cal
S}\overline\otimes\Cal T$ are injective. It then follows easily that
$\Cal S$ (resp\. $\Cal T$) is injective, being a subalgebra of an
injective von Neumann algebra of type II$_1$, and $\Cal M=\Cal
S\overline\otimes\Cal T$ is injective. Since $\Cal S$ and $\Cal T$
appears symmetrically it is enough to prove that $\Cal
S\overline\otimes\Cal R_{\Cal T}$ is injective. In the first place we
will further assume that $\Cal M$ has a finite faithful normal trace
$w$, then we know that we can get a conditional expectation $\gamma$
of norm one from $\Cal M$ onto $\Cal S\overline\otimes\Cal R_{\Cal T}$
 by the formula
 $$\forall\,m\in\Cal M:\gamma(m)\in\Cal S\overline\otimes\Cal R_{\Cal
T}\text{ and }\forall\,x\in\Cal S\overline\otimes\Cal R_{\Cal
T}:w(mx)=w(\gamma(m)x)\;.
 $$

Since there exists a bounded projection $\varphi$ of $B(H)$ onto
$\Cal M$ we get that $\pi=\gamma\circ\varphi$ is a bounded projection
onto $\Cal S\otimes\Cal R_{\Cal T}$. The injective algebra $\Cal
R_{\Cal T}$ contains an increasing sequence of matrix algebras $\Cal
A_k$ such that $\Cal A_0=\Bbb C I$, $\Cal A_k\simeq \Cal M_{2^k}(\Bbb
C)$, $k\geq1$. Letting $\Cal U_0$ denote the unitary group formed as
the union of the unitary groups in $\Cal A_k$ we see that we have an
amenable group and we can mimic the arguments used in the properly
infinite case. Hence we obtain a projection $\pi_0$ onto $\Cal
S\overline\otimes\Cal R_{\Cal T}$ which is a module mapping for the
algebra $\Cal R_0=I_{\Cal S}\otimes(\cup\Cal A_k)$. By assumption
$I_{\Cal S}\otimes\Cal R_{\Cal T}$ is injective and so is the
commutant say $\Cal N$ in $B(H)$, hence there exists a completely
positive projection $\psi$ of $B(H)$ onto $\Cal N$. Since $\pi_0$ is
an $\Cal R_0$-module mapping it follows that $\pi_0(\Cal
N)\subseteq\Cal N$ and  $\pi_0(\Cal N)\subseteq
(\Cal S\overline\otimes\Cal R_{\Cal T})\cap\Cal
N=\Cal S\otimes I_{\Cal R_{\Cal T}}$ on the other hand $\pi_0(\Cal
S\otimes I_{\Cal R_{\Cal T}})=\Cal S\otimes I_{\Cal R_{\Cal T}}$ so
$\pi_0(\Cal N)=\Cal S\otimes I_{\Cal R_{\Cal T}}$. It is now easy to
identify  $(\pi_0\otimes\id):\Cal N\otimes\Cal M_{2^k}(\Bbb C)\to\Cal
S\otimes\Cal M_{2^k}(\Bbb C)$ with $\pi_0:(\Cal N\cup(I_{\Cal S}
\overline\otimes\Cal A_k))''\to\Cal
S\otimes\Cal A_k$. Hence we find that the restriction of $\pi_0$ to
$\Cal N$ is a completely bounded projection of
$\Cal N$ onto $\Cal S\otimes I_{\Cal R_{\Cal T}}$. By defining
$\pi$ on $B(H)$ as $\pi=\pi_0\circ\psi$ we obtain a completely
bounded projection of $B(H)$ onto $\Cal S\otimes I_{\Cal R_{\Cal T}}$
and we have proved that $\Cal S$ is injective. In the same way we get
that $\Cal T$ is injective so $\Cal M\simeq\Cal S\overline\otimes\Cal
T$ is an injective von Neumann  algebra when $\Cal M$ has a finite
faithful normal trace and $\Cal M$ is complemented in $B(H)$. If
$\Cal M$ does not have a faithful finite normal trace we can consider
countably decomposable central projections $e$ in $\Cal S$ and $f$ in
$\Cal T$ and prove by the methods just used that $\Cal M_{e\otimes
f}=\Cal S_e\overline\otimes\Cal T_f$ is injective for each such pair
$(e,f)$ of central projections. The set of such pairs is clearly
upwards filtering and the family $(\Cal M_{e\otimes f})$
of such subalgebras generate $\Cal M$. Hence $\Cal M$ is
injective (may be the easiest way to see this is to remark that the
commutant obviously must have the extension property) and the
corollary follows.$\qquad\square$
\enddemo
\smallskip

\flushpar{\smc Remark 3.3}  Suppose $\Cal M$ is a finite II$_1$ von
Neumann algebra on a Hilbert space $H$ which is a complemented
subspace of $B(H)$. Since any von Neumann subalgebra of such a von
Neumann algebra is complemented, we get that a von Neumann subalgebra
$\Cal S$ of $\Cal M$ is injective if $\Cal S'\cap\Cal M$ is of type
II$_1$.
\vskip.2truecm

The basic content in [4] is well understood and used in [1]. In fact
[1] offers a nice description of a derivation which is inner if and
only if the algebra is injective. This test derivation is basically
the same as introduced in [4] and we will reproduce it once more here
because we think the following presentation is a little bit more
transparent.
\smallskip

\proclaim{Corollary 3.3} Let $\Cal M$ be a von Neumann algebra in
standard position on a Hilbert space $H$ and let $J$ be an antilinear
isometric involution on $H$ such that the commutant $\Cal M'$ equals
$J\Cal M J$. Let $Y$ denote the space of bounded $\Cal M$-module
operators on $B(H)$ which vanishes on $\Cal M$. The space $Y$ is made
into an $\Cal M$ bimodule by the products: $\forall\,m\in\Cal M$,
$\forall\,\varphi\in Y$, $\forall\,x\in B(H)$
 $$(m.\varphi)(x)=\varphi(x)Jm^\ast J\;,\;(\varphi.m)(x)=Jm^\ast
J\varphi(x)\;.
 $$
 Let $\delta:\Cal M\to Y$ be given by $\delta(m)(x)=[Jm^\ast J,x]$,
then $\delta$ is a derivation, and $\Cal M$ is injective if and only
if there exists a $\varphi$ in $Y$ such that
$\delta(m)=\varphi.m-m.\varphi$.
\endproclaim

\demo{Proof} Suppose $\delta=\ad\varphi$ where $\varphi$ is in $Y$ then
 $$\align
[Jm^\ast J,x]&=(\varphi.m-m.\varphi)(x)\\
&=Jm^\ast J\varphi(x)-\varphi(x)Jm^\ast J\\
&=[Jm^\ast J,\varphi(x)]\;,\endalign
 $$
so $[Jm^\ast J,x-\varphi(x)]=0$ for all $x$ in $B(H)$ and all $m$ in
$\Cal M$.
 Hence $\pi(x)$ given by $\pi(x)=x-\varphi(x)$ is a bounded $\Cal M$
module mapping of $B(H)$ onto $\Cal M$, and since $\varphi\mid\Cal
M=0$ we get $\pi\mid\Cal M=\id\mid\Cal M$ so $\pi$ is an $\Cal
M$-module projection of $B(H)$ onto $\Cal M$, and herefore $\Cal M$
is injective by Theorem 3.1. On the other hand if $\Cal M$ is
injective then there exists a completely positive $\Cal M$ module
projection $\pi$ of $B(H)$ onto $\Cal M$, and we can just reverse the
arguments above. Let $\varphi=\id-\pi$ then $\varphi$ is in $Y$ and
$\varphi$ implements the derivation $\delta$.$\qquad\square$

\medskip

{\eightpoint
\flushpar Erik Christensen Matematisk Institut,\newline
Ko/benhavns Universitet,\newline
Universitetsparken 5,\newline
2100 Copenhagen O/\newline
Denmark.}
\medskip

{\eightpoint
\flushpar Allan M\. Sinclair\newline
Department of Mathematics,\newline
University of Edinburgh,\newline
Mayfield Road,\newline
Edinburgh EH9 3JZ,\newline
Scotland.}


\Refs\nofrills{References}

\ref\no[1]\by J\. W\. Bunce, W\. L\. Paschke\paper Quasi Expectations
and Amenable von Neumann Algebras\jour Proc\. Amer\. Math\.
Soc.\vol71\yr1978\pages 232--236
\endref

\ref\no[2]\by E\. Christensen, A.M\. Sinclair\paper On the
Hochschild co\-ho\-mo\-lo\-gy for von Neu\-mann Algebras\jour preprint
\endref

\ref\no[3] \by A\. Connes\paper Classification of injective
factors\jour Ann\. Math\.\vol104\yr1976\pages73--115\endref

\ref\no[4]\by A\. Connes\paper On the Cohomology of Operator
Algebras\jour J\. Funct\. Analysis\vol28\yr1978\pages248--253\endref

\ref\no[5]\by J\. Dixmier\book d'Op\'erateurs dans L'Espace
Hilbertien\publ Gauthier Villars \yr1969\publaddr Paris,
France\endref

\ref\no[6] \by G.A\. Elliott\paper On Approximately Finite Dimensional
von Neumann Algebras II\jour Canad\. Math\. Bull.\vol21\yr1978\pages
415--418\endref

\ref\no[7]\by U\. Haagerup\paper Injectivity and Decomposition of
Completely Bounded Maps, O\-perator Algebras and their Connections with
Topology and Ergodic Theory\jour Springer Lect. Notes in
Math.\vol1132\yr1985\pages170--223
\endref

\ref\no[8]\by R.V\. Kadison\paper Derivations of Operator
Algebras\jour Ann\. Math.\vol83\yr1986\pages280--\newline
293\endref

\ref\no[9]\by R.V\. Kadison, J.R\. Ringrose\book Fundamental of the
Theory of Operator Algebras\publaddr Academic Press\yr1983\publaddr
New York, USA\endref

\ref\no[10] \by V.I\. Paulsen\book Completely bounded Maps and
Dilations\publ Longman Scientific \& Technical\yr1986\publaddr
Haralev, Essex, UK\endref

\ref\no[11]\by G\. Pisier\paper Remarks on Complemented subspaces of
von Neumann Algebras\jour Proc\. Royal Soc.,
Edinburgh\vol121A\yr1992\pages1--4\endref

\ref\no[12] \by G\. Pisier\paper The Operator Hilbert Space and
Complex Interpolation\jour Preprint\endref

\ref\no[13]\by S\. Sakai\paper Derivations of $W^\ast$-algebras\jour
Ann\. Math\.\vol83\yr1966\pages 273--279\endref

\ref\no[14]\by M\. Takesaki\paper Duality in Cross Products and the
Structure of von Neumann Algebras of Type III\jour Acta
Math.\vol131\yr1973\pages 249--310\endref

\ref\no[15] \by J\. Tomiyama\paper On the Projection of Norm One in
$W^\ast$-algebras\jour Proc\. Japan
Acad.\vol33\yr1957\pages608--612\endref

\ref\no[16]\by G\. Wittstock\book Extensions of Completely bounded
$C^\ast$-module Homomor-\newline
phisms\bookinfo  in Proc\. Conference on Operator
Algebras and Group Representations\publ Neptun\yr 1980\moreref
 Pitman\yr1983\publaddr New York, USA\endref
\enddocument